# On Validity of Reed's Conjecture for {$P_5$, Flag$^C$}-free Graphs

Medha Dhurandhar

**Abstract:** Here we prove that Reed's Conjecture is valid for {$P_5$, Flag$^C$ }-free graphs where Flag$^C$ is the complement of the Flag graph. Some of the known results follow as corollaries to our result. Reed's conjecture is still open in general.

**Introduction:**
We consider here simple and undirected graphs. For terms which are not defined herein we refer to Bondy and Murty [1]. In 1998, Reed proposed the following Conjecture which gives, for any graph G, an upper bound for its chromatic number χ(G) in terms of the clique number ω(G) and the maximum degree Δ(G).

**Reed's Conjecture [2]:** For any graph G, $\chi(G) \leq \lceil \frac{\Delta + \omega + 1}{2} \rceil$.

In [3], [4], [5], [6], [7], [8], [9], [10], [11], [12], [13], [14] it is shown that Reed's Conjecture holds for some graph classes defined by forbidden configurations:
• ($P_5$, $P_2 \cup P_3$, House, Dart)-free graphs,
• ($P_5$, Kite, Bull, ($K_3 \cup K_1$) + $K_1$)-free graphs,
• ($P_5$, C4)-free graphs,
• (Chair, House, Bull, $K_1 + C_4$)-free graphs,
• (Chair, House, Bull, Dart)-free graphs.
• $3K_1$-free graphs
• {$2K_2$, $C_4$}-free graphs
• Quasiline graphs
• $K_{1,3}$-free
• Generalized line graphs
• Graphs with $\chi \leq \omega + 2$
• Planar and toroidal graphs
• Decomposable graphs
• Perfect graphs
• Line graphs of Multigraphs
• Graphs with disconnected complements
• Graphs G with $\chi(G) > \lceil \frac{V(G)}{2} \rceil$ and graphs G with $\Delta(G) > \lceil \frac{V(G) - \alpha(G) + 3}{2} \rceil$
• Graphs G with $\Delta(G) \geq |V(G)| - 7$, and graphs G with $\Delta(G) \geq |V(G)| - \alpha(G) - 4$

**Notation:** For a graph G, V(G), E(G), Δ(G), ω(G), χ(G) denote the vertex set, edge set, maximum degree, size of a maximum clique, chromatic number respectively. For u ∈ V(G), N(u) = {v ∈ V(G) / uv ∈ E(G)}, and $\overline{N(u)}$ = N(u)∪(u). If S ⊆ V(G), then <S> denotes the subgraph of G induced by S. Also for u ∈ V(G), deg u denotes the degree of u in G. If C is some coloring of G and if a vertex u of G is colored m in C, then u is called a m-vertex. All graphs considered henceforth are simple.

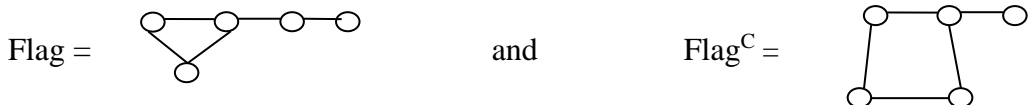

Flag =   and   Flag$^C$ =

**Figure 1**

This paper proves that Reed's Conjecture holds for {$P_5$, Flag$^C$}-free graphs.

**Note:** As G is $P_5$-free, the only odd, chordless cycles in G are $C_5$.

**Theorem:** If G is a {$P_5$, Flag$^C$}-free, then $\chi(G) \leq \lceil \frac{\Delta + \omega + 1}{2} \rceil$.

Proof: Let G be a smallest $\{P_5, \text{Flag}^C\}$-free graph with $\chi(G) > \lceil \frac{\Delta + \omega + 1}{2} \rceil$. Let $u \in V(G)$. Then by minimality, $\chi(G) - 1 \leq \chi(G-u) \leq \lceil \frac{\Delta(G-u) + \omega(G-u) + 1}{2} \rceil \leq \lceil \frac{\Delta + \omega + 1}{2} \rceil < \chi(G)$. Thus $\chi(G-u) = \chi(G) - 1 = \lceil \frac{\Delta + \omega + 1}{2} \rceil \ \forall \ u \in V(G)$. Let C be a $\chi(G)$-coloring of G and R = $\{x \in N(u)/ \ x$ receives a unique color from C in $<N(u)>\}$. Then $\Delta \geq \deg u \geq |R| + 2(\lceil \frac{\Delta + \omega + 1}{2} \rceil - |R|)$ and $|R| \geq \omega(G) + 1$. **I**

Further let R = S∪T where S = $\{x \in R/ \ xy \in E(G) \ \forall \ y \in R-x\}$ and T = R-S. Also let T' = $\{x' \in V(G) - \overline{N(u)} / \exists \ x, y \in T$ with $xy \notin E(G), \ x'x \in E(G)$ and $x', y$ have the same color$\}$. Further let $S_1$ = $\{w \in S-T/ \exists \ y \in T'$ with $wy \notin E(G)\}$ and $S_1$' = $\{x \in V(G) - \overline{N(u)} / $ for some $\alpha$-vertex $y \in S_1$ and $z \in T'$ s.t. $yz \notin E(G), x$ is an $\alpha$-vertex of $z\}$.

We have

1. $|T| \geq 2$ (else $|T| = 0$ and by **I**, $<R> \supseteq K_{\omega(G)+1}$) and $|S_1| > 0 \Rightarrow |S_1'| > 0$.

2. **If $t, t' \in T$ s.t. $tt' \notin E(G)$ and $t, t'$ have colors $j, i$ resply, then $t$ ($t'$) has a unique $i$-vertex ($j$-vertex).**

Since $t$ ($t'$) is a unique $j$ ($i$)-vertex of $u$, $\exists$ a $j$-$i$ path say P = $\{t, V, W, t'\}$. Let if possible $t$ have another $i$-vertex V'. Then if V'W $\notin E(G)$, $<W, t', u, t, V'> = P_5$, and if V'W $\in E(G)$, $<V, W V', t, u> = \text{Flag}^C$, a contradiction in both the cases.

3. **If $t, t', t'' \in T$ s.t. $tt', tt'' \notin E(G)$, and $t$ has color $i$, then $t', t''$ have a common, unique $i$-vertex.**

Let $t', t''$ have colors $j, k$ resply. By **2**, each of $t', t''$ has a unique $i$-vertex say A, B resply. Let if possible $A \neq B$ and P = $\{t, V, A, t'\}$, M = $\{t, W, B, t''\}$ be the $i$-$j$, $i$-$k$ paths. Then $t't'' \in E(G)$ (else $<A, t', u, t'', B> = P_5$). But then if $t'W \in E(G)$, $<t', t'', B, W, t> = \text{Flag}^C$ and if $t'W \notin E(G)$ $<t, W, B, t'', t'> = P_5$, a contradiction in both the cases.

4. **$<S_1' \cup T'>$ is complete.**

**First we prove that $<T'>$ is complete.** Let if possible $\exists \ a', b' \in T'$ s.t. $a'b' \notin E(G)$. Then by definition of T' $\exists \ A, B, C, D \in T$ s.t. A, a'; B, b' have same colors and AC, BD $\notin E(G)$. Clearly AB $\in E(G)$ (else by **3**, Ba', Ab' $\in E(G)$ and $<b', A, u, B, a'> = P_5$). Also AD $\in E(G)$ (else a'd', a'D $\in E(G)$ and $<a', d', b', D, u> = \text{Flag}^C$). Similarly BC $\in E(G)$. Now a'B $\in E(G)$ (else a'D $\in E(G)$ as otherwise $<D, A, B, C, a'> = \text{Flag}^C$ or $P_5$ depending on whether CD $\in E(G)$ or $\notin E(G)$. Also $<a', d', b', D, u> = \text{Flag}^C$ if a'd' $\in E(G)$ and $<B, d', b', D, a'> = P_5$ if a'd' $\notin E(G)) \Rightarrow$ a'D $\notin E(G)$ (else $<u, B, a', D, b'> = \text{Flag}^C$). But then $<a', B, u, D, b'> = P_5$, a contradiction.

**Next we prove that $st' \in E(G) \ \forall \ s \in S_1', \ t' \in T'$.** Let if possible $\exists \ s' \in S_1'$ and $a' \in T'$ s.t. $s'a' \notin E(G)$. By **I**, $\exists \ A, B \in T$ s.t. AB $\notin E(G)$ and $\{A, b', a', B\}$ is a bi-color path. Also s'b' $\notin E(G)$ (else $<s', b', a', B, u> = \text{Flag}^C$ or $P_5$ according as whether s'B $\in E(G)$ or $\notin E(G)$). Now by definition let $s \in S$ be an $\alpha$-vertex and c' $\in T'$ be s.t. sc' $\notin E(G)$ and s' is the $\alpha$-vertex of c'. As $<T'>$ is complete, ca', cb' $\in E(G)$. Also both s'B and s'A $\notin E(G)$ (else if s'A, s'B $\in E(G)$ then $<s', A, u, B, a'> = \text{Flag}^C$ and if only s'A $\in E(G)$, then $<a', B, u, A, s'> = P_5$). Now sa' $\notin E(G)$ (else $<u, s, a', c', s'> = P_5$). Similarly sb' $\notin E(G) \Rightarrow$ c'B $\in E(G)$ (else $<s, B, a', c', s'> = P_5$) and similarly c'A $\in E(G)$. But then $<B, u, A, c', s'> = \text{Flag}^C$, a contradiction.

**Finally we prove that $<S_1'>$ is complete.** Let if possible $\exists \ s_{11}', s_{12}' \in S_1'$ s.t. $s_{11}'s_{12}' \notin E(G)$. Now by **1**, $\exists \ A, B \in T; a', b' \in T'$ s.t. AB $\notin E(G)$ and $\{A, b', a', B\}$ is a bi-color path. As $st' \in E(G) \ \forall \ s \in S_1', t' \in T'$, $s_{11}'a', s_{12}'a', s_{11}'b', s_{12}'b' \in E(G)$. Further either $s_{11}'B \notin E(G)$ or $s_{12}'B \notin E(G)$ (else $< s_{11}', b', s_{12}', B, u> = \text{Flag}^C$). W.l.g. let $s_{11}'B \notin E(G)$. Then $s_{11}'A \in E(G)$ (else $<B, u, A, b', s_{11}'> = P_5$) $\Rightarrow s_{12}'A \notin E(G)$ (else $<s_{11}', a', s_{12}', A, u> = \text{Flag}^C$) and $s_{12}'B \in E(G)$ (else $<B, u, A, a', s_{12}'> = P_5$). But then $< s_{11}', A, u, B, s_{12}'> = P_5$, a contradiction.

Let $T = S_0$ and $T' = S_0'$. Further define $S_l = \{x \in S - \bigcup_{m=0}^{l-1} S_m / \exists\, y \in S_{l-1}'$ with $xy \notin E(G)\}$ and $S_l' = \{z \in V(G) - \overline{N(u)} /$ for some $\alpha$-vertex $x \in S_l$ and $y \in S_{l-1}'$ s.t. $xy \notin E(G)$, $z$ is an $\alpha$-vertex of $y\}$, $1 \leq l$. By construction every color used in $S_l$ is used in $S_l'$, $l \geq 0$. Let $S_0', .., S_k'$ be a maximal such sequence. Then clearly every $x \in W = S - \bigcup_{l=0}^{k-1} S_l$ is adjacent to every vertex of $\bigcup_{l=0}^{k} S_l'$.   **II**

**Claim:** $\langle W \cup \{\bigcup_{l=0}^{k} S_l'\} \rangle$ is complete.

We prove the **Claim** by induction. By **4**, the result is true for $i = 0, 1$. Assume that $\langle \bigcup_{l=0}^{k-1} S_l' \rangle$ is complete. Let if possible $\exists\, x \in S_k'$, $y \in S_j'$ s.t. $xy \notin E(G)$. Let $j$ be the smallest such index. By construction $\exists$ $\beta$-vertex $z \in S_k$ and $w \in S_{k-1}'$ s.t. $zw \notin E(G)$, $x$ is the $\beta$-vertex of $w$. Now $j > 0$ (else by induction $wy \in E(G)$, by construction $zy \in E(G)$ and $\langle u, z, y, w, x \rangle = P_5$). Let $a' \in S_0'$. Then by definition $\exists$ $A, B \in T$ and $b' \in S_0'$ s.t. $\{A, b', a', B\}$ is a bi-color path. As $j > 0$, $xa', xb' \in E(G)$. Also by induction $ya', yb' \in E(G)$. Now both of $xA, yA$ ($xB, yB$) $\notin E(G)$ (else $\langle x, a', y, A, u \rangle = \text{Flag}^C$). But $xA$ or $xB$ ($yA$ or $yB$) $\in E(G)$ (else $\langle B, u, A, b', x \rangle = P_5$). W.l.g. let $xA \in E(G) \Rightarrow yA \notin E(G)$ (else $\langle x, a', y, A, u \rangle = \text{Flag}^C$). Then $yB \in E(G)$ (else $\langle B, u, A, b', y \rangle = P_5$) $\Rightarrow xB \notin E(G)$ (else $\langle x, b', y, B, u \rangle = \text{Flag}^C$). But then $\langle x, A, u, B, y \rangle = P_5$, a contradiction. Thus $\langle \bigcup_{l=0}^{k} S_l' \rangle$ is complete.

Hence by **II**, the **Claim** holds.

But then as every unique color used in $\langle N(u) \rangle$ is used in $\langle W \cup \{\bigcup_{l=0}^{k} S_l'\} \rangle$, by **I**, $\langle W \cup \{\bigcup_{l=0}^{k} S_l'\} \rangle \supseteq K_{\omega(G)+1}$, a contradiction.

This proves the Theorem.

**Corollary 1**: If G is $\{P_5, C_4\}$-free, then $\chi(G) \leq \lceil \frac{\Delta + \omega + 1}{2} \rceil$.

**Corollary 2**: If G is $3K_1$-free, then $\chi(G) \leq \lceil \frac{\Delta + \omega + 1}{2} \rceil$.

**Corollary 3**: If G is $\{P_3 \cup K_1\}$-free, then $\chi(G) \leq \lceil \frac{\Delta + \omega + 1}{2} \rceil$.

**Corollary 4**: If G is $\{2K_2, C_4\}$-free, then $\chi(G) \leq \lceil \frac{\Delta + \omega + 1}{2} \rceil$.

**References**

[1] J.A. Bondy and U.S.R. Murty. Graph Theory, volume 244 of Graduate Text in Mathematics, Springer, 2008.

[2] B. Reed, "$\omega, \Delta$, and $\chi$", Journal of Graph Theory, 27:177–212, 1998.

[3] N.R. Aravind, T. Karthick, and C.R. Subramanian, "Bounding $\chi$ in terms of $\Delta$ and $\omega$ for some classes of graphs", Discrete Mathematics, 311:911–920, 2011.

[4] M. Chudnovsky, N. Robertson, P.D. Seymour, and R. Thomas, "The Strong Perfect Graph Theorem", Annals of Math., 164:51–229, 2006.

[5] J.L. Fouquet, V. Giakoumakis, F. Maire, and H. Thuillier, "On graphs without P5 and P5", Discrete Mathematics, 146:33-44, 1995.

[6] J.L. Fouquet and J.M. Vanherpe, "Reed's conjecture on hole expansions", Technical report, L.I.F.O., 2011.

[7] L. Rabern, "A note on Reed's conjecture", SIAM Journal on Discrete Mathematics, 22:820–827, 2008.

[8] D. Seinsche, "On a property of the class of n-colorable graphs", Jour. of Comb. Theory, Series B (16):191–193, 1974.

[9] Ingo Schiermeyer, "Chromatic number of $P_5$-free graphs: Reed's conjecture", Discrete Mathematics, 2015.

[10] L. Esperet, L. Lemoine, F. Maffray, G. Morel, "The chromatic number of {P5, K4}-free graphs", Discrete Math. 313:743–754, 2013.

[11] A. Kohl, I. Schiermeyer, "Some results on Reed's conjecture about ω,Δ and χ with respect to α", Discrete Math. 310:1429–1438, 2010.

[12] L. Rabern, "A note on Reed's conjecture", SIAM J. Discrete Math. 22 (2): 820–827, 2008.

[13] A.D. King, B. Reed, A. Vetta, "An upper bound for the chromatic number of line graphs", European J. Comb. 28: 2182–2187, 2007.

[14] A. Kohl, I. Schiermeyer, Some results on Reed's conjecture about ω,Δ and χ with respect to α, Discrete Math. 310, 1429–1438, 2010.